\documentclass[11pt,a4paper]{article}

\usepackage{graphicx}

\usepackage[utf8]{inputenc}
\usepackage{dsfont}
\usepackage{textgreek,textcomp}
\usepackage{amsmath,bm}
\usepackage{amssymb}
\usepackage{algorithm}
\usepackage{stmaryrd}
\usepackage{geometry}
\usepackage[noend]{algpseudocode}
\usepackage{float,graphicx,tikz,tikz-3dplot,pgfplots}
\usetikzlibrary{calc, shapes.geometric, arrows, angles, quotes}
\tikzstyle{arrow} = [thick,->,>=stealth,-latex]
\usepackage{graphicx}
\usepackage{subcaption}
\usepackage{xcolor}
\usepackage[normalem]{ulem}

\pgfplotsset{compat=1.16}

\begin{document}
	
	\title{Advanced Modeling of Rectangular Waveguide {Devices with Smooth Profiles} by Hierarchical Model Reduction}

	\author{{Gines} Garcia-Contreras$^\dagger$, {Juan} C\'orcoles$^\dagger$, {Jorge A.} Ruiz-Cruz$^\dagger$, \\
		{Matteo} Oldoni$^\diamond$, {Gian Guido} Gentili$^\diamond$, {Stefano} Micheletti$^\#$, {Simona} Perotto$^\#$}
	\maketitle 
	\noindent
		{\small $^\dagger$ Departamento de Tecnolog\'ia Electr\'onica y de las Comunicaciones, Escuela Polit\'ecnica Superior, Universidad Aut\'onoma de Madrid, 28049 Madrid, Spain\\
			 {\tt \{gines.garcia,juan.corcoles,jorge.ruizcruz\}@uam.es}\\
			$^\diamond$ Dipartimento di Elettronica, Informazione e Bioingegneria, 
			 Politecnico di Milano, Piazza L. da Vinci 32, I-20133 Milano, Italy
			{\tt \{gianguido.gentili,matteo.oldoni\}@polimi.it}\\[2mm]
			$^\#$ MOX - Dipartimento di Matematica,
			Politecnico di Milano, Piazza L. da Vinci 32, I-20133 Milano, Italy
			{\tt \{stefano.micheletti,simona.perotto\}@polimi.it}\\[2mm]
		}
	
	\date{}
	
	\maketitle
	
\begin{abstract}
	We {present a new method for the analysis of smoothly varying tapers, transitions and filters 
		in rectangular waveguides. With this aim, we apply a
		Hierarchical Model (HiMod) reduction 
		to the vector Helmholtz equation. We exploit a suitable coordinate transformation and, successively, we use the waveguide modes as a basis for the HiMod expansion. 
		We show that accurate results can be obtained with an impressive speed-up factor when compared with standard commercial codes based on a three-dimensional
		finite element discretization.}
\end{abstract}

\textbf{Keyords} Waveguide devices, Modal analysis, Hierarchical Model reduction, Finite Element Method.

%

\section{Introduction}
%
%
%
%
{H}{ierarchical} model (HiMod) reduction is a powerful numerical tool to speed up the analysis of complex structures in the presence of directional features. HiMod reduction has been successfully employed in a variety of different application fields, such as in linear acoustics \cite{GP}, for advection-diffusion-reaction phenomena \cite{HM1,HM2,Perotto,SP14,LupoSimo22,LupoSimo23}, for the blood flow modeling~\cite{GuzzettiPerotto18}, also in patient-specific artery segments \cite{Perotto2}. Unlike typical model order reduction techniques, such as \cite{Wu, Kulas}, 
the underlying key advantage of the HiMod approach is that the original problem size is shrunk through the application of a different discretization along the main and the secondary direction, according to a separation of variable principle. This strategy allows accurately representing the problem dynamics without introducing a heavy model simplification. In particular, a careful selection of the bases supporting the different discretizations allows one for an extensive saving on computer resources, both in terms of processing time and memory storage. The 
main challenge of a HiMod reduction is, therefore, being able to find tailored bases to discretize the problem along the two directions. 

Propagation systems where the structure extends mainly along a single direction of space, such as waveguide tapers, smooth transitions, twists and bends \cite{uher1993waveguide}, represent an ideal environment for the application of a HiMod approach. In particular, modes offer a straightforward solution to the problem of finding a suitable reduced basis for the modeling of transverse dynamics, thanks to the completeness of modal expansion in the transverse plane, perpendicular to the propagation direction \cite{collin2007foundations}. Moreover, a modal representation turns out to be extremely effective given that the transverse field can be well approximated by a reduced amount of modes. If such a condition holds, the modal spectrum is a natural and extremely attractive reduced basis for a HiMod approximation of the transverse dynamics. Nevertheless, we observe that the direct application of the modal expansion might be unpractical since the profile of the device in the transverse plane may change,
the modal fields becoming a function of the shape itself. 

The method proposed in this work overcomes this limitation by introducing a coordinate transformation that converts the actual structure into a geometry with a constant shape in the transverse plane, invariant along the propagation direction. This expedient allows 
resorting to a \textit{unique} set of modes to represent the transverse features of the field along the whole domain.
Of course, in order to make this representation equivalent to the original problem, the material inside the waveguide has to be substituted with a new artificial medium, which is directly obtained from the coordinate transformation. This approach has often been referred to as \textit{Transformation optics} or \textit{Transformation electromagnetics}, having several applications in the modelling and design of passive components, antennas, metasurfaces and cloaking \cite{Chen,Sun,Kuang,Pendry}.
\begin{figure}[!t]
	\centering
	\subfloat[]{\includegraphics[width=0.47\columnwidth]{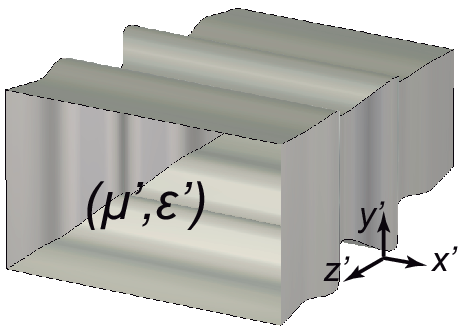}%
		\label{subfig:t1}}
	\subfloat[]{\includegraphics[width=0.47\columnwidth]{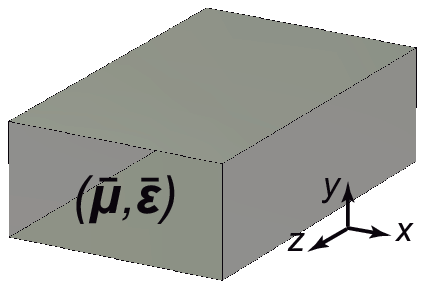}%
		\label{sufig:t2}}
	\caption{Geometric transformation. (a) Original domain with scalar permittivity and permeability. The waveguide exhibits a mild width and height variation along the $z'$-axis. (b) Transformed domain with tensorial permittivity and permeability. The profile is now constant along the $z$-axis. }
	\label{fig:transformation}
\end{figure}

The HiMod technique here proposed is applied to tapers, smooth transitions and filters in rectangular waveguides. The main focus of tapered structures is illustrated by wide-band matched transitions between dissimilar rectangular waveguides
\cite{Sangam,Garner}. Other common applications are found in mode converters \cite{Liu,Ren} and pyramidal horns \cite{Rebo,Otter}.

Classic papers present design and analysis procedures based on solving the three-dimensional (3D) problem associated with the taper by using Mode-Matching (MM) \cite{Rebo}. More recently, the Finite Element Method (FEM) was applied through the use of well-known commercial software \cite{Ren}. For smooth walled tapers, the extremely efficient MM technique can be less engaging due to the need of representing the geometry by a large number of steps, thereby increasing the computational time and introducing an undesired geometry discretization and approximation. Thus, FEM, in the 3D version, has become until today the preferred choice, even though the computing time can become an issue when dealing with long tapers.  

The method proposed in this work offers a totally new concept, since it combines the powerful, well-known, modal expansion to describe the transverse dynamics with a standard one-dimensional (1D) FE discretization in order to model the phenomenon along the propagation direction. The extremely efficient representation of the problem makes HiMod reduction the best candidate for ``long'' structures (with regards to the wavelength) with smoothly varying cross-sections. The examples of application discussed in this paper provide just a glance at the great potentiality that this technique can provide.

The theoretical aspects of the problem and, in particular, the coordinate transformation are discussed in Section \ref{sec2}. In Section \ref{sec3} we introduce 
the HiMod reduction. Some issues regarding the function spaces are discussed, together with a detailed analysis of the effect of the coordinate transformation on the output port and the corresponding modal functions. The final purpose of this analysis is, in fact, the computation of the multimode Generalized Scattering Matrix (GSM) \cite{guglielmi1999advanced}. This is achieved by introducing some suitably normalized mode functions as forcing terms in the HiMod model. Section \ref{sec4} presents the results obtained by applying the new approach {to tapered structures}, while in the final section we draw some conclusions.

\section{Problem statement}\label{sec2}
We discuss the application of a HiMod reduction to the structure represented in Fig.~\ref{fig:transformation}(a). 
The process can be split into two phases: To begin with, we carry out a coordinate transformation so as to convert the tapered geometry into a uniform section waveguide of length $L$. In completion of this step, a corresponding artificial material is introduced to make the electromagnetic problem equivalent to the original one with tapered geometry. In the second step, we carry out a discretization that employs the standard modes of the uniform structure as basis functions in the transverse plane, plus a classic 1D-FEM discretization along the propagation direction. 
This section focuses on the first step.


We consider a rectangular waveguide device in the $(x',y',z')$ reference frame so that the cross-section, lying in the $x'y'$ plane, exhibits smooth variations in the profile (see Fig.~\ref{fig:transformation}(a)). The taper is oriented such that port $1$, of size $a_0$ and $b_0$, along the $x'$- and $y'$- direction, respectively, lies at $z'=0$, while port $2$, of size $a_L$ and $b_L$ along the $x'$- and $y'$-direction, respectively, lies at $z'=L$. The profile of the device is, then, described by two functions, $a = a(z')$ and $b = b(z')$, depending on the propagation direction, $z'$, so that 
\begin{equation}\label{bc_ab}
a(0) = a_0, \ b(0) = b_0, \ a(L) = a_L, \ b(L) = b_L.
\end{equation}  
Then, we define a geometric transformation 
\begin{equation}\label{rrp}
{\bf r} = {\bf r({\bf r}')},
\end{equation}
with ${\bf r}' = (x',y',z')$ and ${\bf r}=(x,y,z)$, such that the original device becomes a waveguide of constant rectangular section, $a_0 \times b_0$, for all $z'$, as illustrated in Fig. \ref{fig:transformation}(b). 
It is assumed that port $1$ is defined by $-a_0/2 < x' < a_0/2$ and $-b_0/2 < y' < b_0/2$, while port $2$ is characterized by the ranges $-a_L/2 < x' < a_L/2$ and $-b_L/2 < y' < b_L/2$. The taper geometry functions $a=a(z')$, $b=b(z')$, can be quite arbitrary, the only requirements being the endpoint conditions in \eqref{bc_ab}. Here, we select the transformation \eqref{rrp} as follows:
\begin{equation}
\left\{
\begin{array}{lll}
\displaystyle x = x'\frac{a_0}{a(z')}, & -a(z')/2\le x' \le a(z')/2, & 0 \leq z' \leq L\\[4mm]
\displaystyle y = y'\frac{b_0}{b(z')}, & -b(z')/2\le y' \le b(z')/2,  & 0 \leq z' \leq L\\[4mm]
z = z', & & 0 \leq z' \leq L.
\end{array}
\right.  \label{eq:transformation}
\end{equation}
It can be easily checked that in the new frame $(x,y,z)$, the taper geometry has now a constant section independent of $z$, i.e., the physical domain has been transformed into the rectangular prism $\Omega$ of dimensions $a_0 \times b_0 \times L$. 

In order to ensure that the electromagnetic problem in the original and in the transformed geometry are equivalent, we have to define a non-isotropic, non-homogeneous artificial material filling the waveguide with the following properties:
\begin{align}
&\bar{\bm{\epsilon}}_r = {\epsilon}'_r\frac{\bar{\bf{J}}\;\bar{\bf{J}}^T}{\det(\bar{\bf{J}})}=\epsilon'_r\bar{\bm\Lambda},\label{permit}\\
&\bar{\bm{\mu}}_r =  {\mu}_r'\frac{{\bar{\bf{J}}}\;{\bar{\bf{J}}}^T}{\det(\bar{\bf{J}})} = \mu'_r\bar{\bm\Lambda},\label{permea}
\end{align}
where $\epsilon'_r$, $\mu'_r$ are the homogeneous relative permittivity and permeability in the original reference frame, respectively, 
\begin{equation}
\bar{\bf{J}} =\frac{\partial{\bf r}}{\partial{\bf r}'}
=
\begin{bmatrix}
\frac{\partial x}{\partial x'} & \frac{\partial x}{\partial y'} & \frac{\partial x}{\partial z'} \\[2mm]
\frac{\partial y}{\partial x'} & \frac{\partial y}{\partial y'} & \frac{\partial y}{\partial z'} \\[2mm]
\frac{\partial z}{\partial x'} & \frac{\partial z}{\partial y'} & \frac{\partial z}{\partial z'} 		
\end{bmatrix}, \label{eq:Jacobian}
\end{equation}
is the Jacobian of the transformation, {and 
	$\bar{\bm\Lambda}=\bar{\bf{J}}\;\bar{\bf{J}}^T/\det(\bar{\bf{J}})$ denotes the anisotropy  contribution in the material originated from the coordinate transformation}.
By substituting \eqref{eq:transformation} into \eqref{eq:Jacobian}, we obtain:
\begin{equation}
\bar{\bf{J}} = \begin{bmatrix}
\frac{a_0}{a(z')} & 0 & -\frac{x}{a(z')}\frac{da(z')}{dz'} \\[2mm]
0 & \frac{b_0}{b(z')} & -\frac{y}{b(z')}\frac{db(z')}{dz'} \\[2mm]
0 & 0 & 1 		
\end{bmatrix}.
\end{equation}
Solving the equivalent electromagnetic problem in the uniform waveguide section with non-homogeneous anisotropic permittivity and permeability yields the transformed electric and magnetic field, ${\bf E}$ and ${\bf H}$, respectively which are related to the corresponding fields in the original problem, ${\bf E}'$, ${\bf H'}$, by the relations
\begin{equation}\label{eqEH}
{\bf{E}} = \bar{\bf{J}}^T{\bf{E}}', \quad
{\bf{H}} = \bar{\bf{J}}^T{\bf{H}}'. 
\end{equation}

The technique described above is quite popular and can lead to great computational benefits in the analysis of complex structures \cite{GP}. 
In our case, the coordinate transformation \eqref{eq:transformation} represents only the first step in building the reduced model we are interested in, as discussed in the following section.

\section{{HiMod} Discretization}\label{sec3}

\subsection{The reduced basis}
\begin{figure}[!t]
	\centering
	\subfloat[]{\includegraphics[width=0.35\columnwidth]{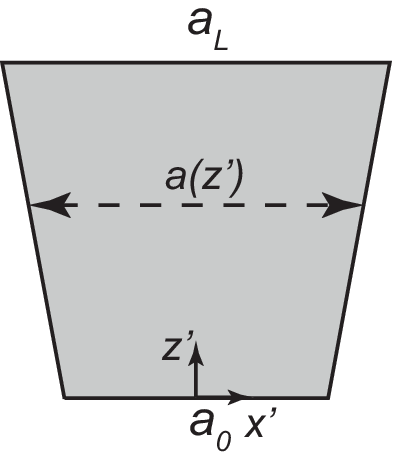}%
		\label{subfig:1}}
	\subfloat[]{\includegraphics[width=0.43\columnwidth]{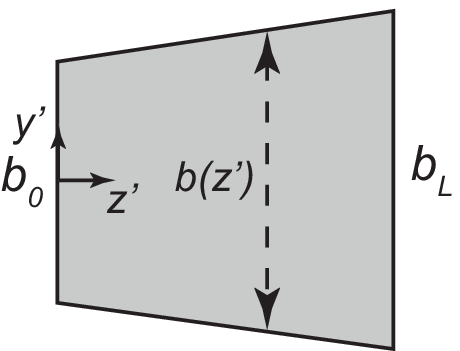}%
		\label{sufig:2}} \\
	\subfloat[]{\includegraphics[width=0.38\columnwidth]{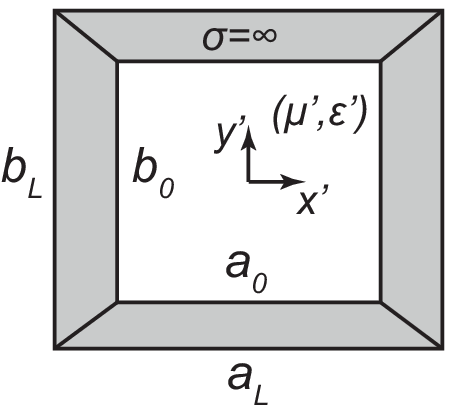}%
		\label{subfig:3}}
	\hfil
	\subfloat[]{\includegraphics[width=0.35\columnwidth]{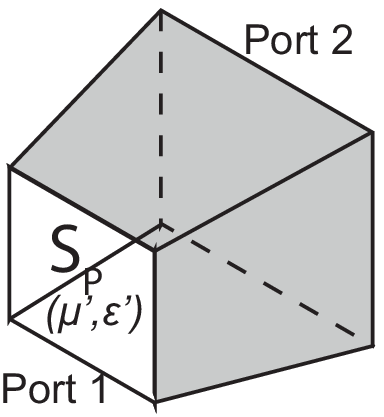}%
		\label{subfig:4}}
	\caption{Generic taper in the original domain. (a) Top view, (b) right side view, (c) front view, (d) isometric view.}
	\label{fig:perspective}
\end{figure}

With a view to the application of {HiMod}, a convenient  basis in the {transformed} geometry for the representation of the electromagnetic field in the transverse plane is represented by the set of Transverse Electric (TE) and Transverse Magnetic (TM) modes, describing the electric field in a uniform and homogeneously-filled rectangular waveguide. 
{With reference to the configuration in Fig.~\ref{fig:perspective}, the modal
	basis will be} uniquely defined from the knowledge of the geometry of port $1$, and it can be used to represent the electric field anywhere in the waveguide profile thanks to the geometric transformation in {\eqref{eq:transformation}}. 

We highlight that modal functions exactly represent the modes at port $1$, since, by definition, no shape modification takes place in the plane where the port lies. Thus, for the same reason, such modal functions do not represent, in general, the modes at port $2$, nor at any internal cross-section of the transformed geometry. More importantly, 
TE and TM modes
provide a complete basis for any field satisfying Maxwell's equations, including the non-homogeneous equivalent waveguide at hand. {Thus, TE and TM modal functions provide a natural and convenient reduced basis for a HiMod approximation.}

By splitting the electric field into a transverse (subscript $t$) and a longitudinal (subscript $z$) component, we {can rewrite ${\bf E}$ as ${\bf E}(x,y,z)={\bf{E}}_t(x,y,z)+{\bf E}_z(x,y,z)$, with}
\begin{align}
&{\bf{E}}_t(x,y,z) = \sum_{n=1}^{N_{\rm TE}}{\bf{e}}^{\rm TE}_{t,n}(x,y)\tau^{\rm TE}_{n}(z) 
+\sum_{n=1}^{N_{\rm TM}}{\bf{e}}^{\rm TM}_{t,n}(x,y)\tau^{\rm TM}_{n}(z) \label{eq:et0} \\
&{\bf E}_z(x,y,z) = \sum_{n=1}^{N_{\rm TM}}{\bf e}^{\rm {TM}}_{z,n}(x,y)\zeta_{n}(z), \label{eq:ez0}
\end{align}
where ${\bf e}^{\rm TE}_{t,n}$, ${\bf e}^{\rm TM}_{t,n}$ are transverse modal fields, ${\bf e}^{\rm TM}_{z,n}$ are longitudinal modal fields, while functions $\tau^{\rm TE}_{n}(z)$, $\tau^{\rm TM}_{n}(z)$ and $\zeta_{n}(z)$ model the $z$-dependence of the electric field components. 
All modal fields ${\bf e}^{\rm TE}_{t,n}$, ${\bf e}^{\rm TM}_{t,n}$, ${\bf e}^{\rm TM}_{z,n}$ can be chosen as real functions and can be obtained from the usual mode decomposition of a rectangular waveguide, including the perfect conductor boundary condition assigned on the lateral surface of the taper~\cite{Pozar}.
Their superposition as functions of $z$ is instead determined by $\tau^{\rm TE}_{n}(z)$, $\tau^{\rm TM}_{n}(z)$ and $\zeta_{n}(z)$ functions, still unknown and objective of the solution process.  Note that $\tau_n$ and $\zeta_n$ are not necessarily equal due to the artificial material inside the transformed domain.
Expressions \eqref{eq:et0}-\eqref{eq:ez0} represent the electric field at all points in the transformed geometry. The associated theoretically infinite series have been truncated to $N_{\rm TE}$ and $N_{\rm TM}$ terms for the ${\rm TE}$ and ${\rm TM}$ modes, respectively. It will be shown that, in practice, a relatively small value of such terms is sufficient to obtain an accurate representation of field ${\bf{E}}$ in the structure. Since modes in the rectangular waveguide are identified by two indexes, $p$ and $q$, associated with $x$ and $y$, respectively (see the Appendix), one has to define the correspondence of the generic index $n$ to the indexes $p$, $q$, namely 
\begin{equation}\label{pnindices}
p=p(n)=p_n, \quad q=q(n)=q_n.
\end{equation}

The $z$-dependence of the fields in \eqref{eq:et0}-\eqref{eq:ez0} can be represented in several ways. 
{We are interested in a finite element discretization, due to the relevant flexibility in terms of polynomial degree, to the high accuracy and to the remarkable efficiency guaranteed by the sparsity pattern characterizing the associated matrices. For this reason, we select
	the 1D interpolating Lagrange polynomials to expand the function in \eqref{eq:et0}-\eqref{eq:ez0} depending on $z$, so that we have}
\begin{align}\label{eqtetau}
\tau_{n}^{\rm TE}(z) &= \sum_{l=1}^{N_{l,t}}c^{\rm TE}_{n,l} \phi_{l}(z) \\
\label{eqtmtau}
\tau_{n}^{\rm TM}(z) &= \sum_{l=1}^{N_{l,t}}c^{\rm TM}_{n,l} \phi_{l}(z) \\
\label{eqtmz}
\zeta_{n}(z) &= \sum_{l=1}^{N_{l,z}}d^{\rm TM}_{n,l} \psi_{l}(z),
\end{align}
where $\{\phi_{l}(z)\}$ and $\{\psi_{l}(z)\}$ are the standard Lagrange polynomials providing a basis of {dimension $N_{l,t}$ and $N_{l,z}$, respectively} for a 1D FE discretization (we use different symbols for the transverse and longitudinal polynomials for reasons that will be clear below), whereas $c_{n,l}^{\rm TE}$, $c_{n,l}^{\rm TM}$ and $d_{n,l}^{\rm TM}$ denote the corresponding degrees of freedom (DOFs), i.e., the actual unknowns. From \eqref{eq:et0}-\eqref{eq:ez0} and \eqref{eqtetau}-\eqref{eqtmz}, it follows that the total number of DOFs used to represent the electric field in the structure is equal to
\begin{equation*}
N_{tot} = (N_{\rm TE}+N_{\rm TM})N_{l,t}+N_{\rm TM}N_{l,z}.
\end{equation*}

In Section~\ref{sec4}, we will verify 
that the computational advantage provided by a HiMod reduction is the much smaller global number of DOFs when commpared  with a full 3D FE discretization.

\subsection{Weak form and {algebraic HiMod formulation}}
{The model we reduce does coincide with Maxwell's equations, which,}
in the transformed domain $(x,y,z)$, read as:
\begin{align}
&{\nabla} \times {\bf{E}} = -\jmath\omega \mu_0\bar{\bm\mu}_r{\bf{H}}, \label{eq:maxwell1}\\
&{\nabla} \times {\bf{H}} = \jmath\omega \epsilon_0\bar{\bm\epsilon}_r{\bf{E}}, \label{eq:maxwell2}
\end{align}
where $\omega$ is the angular frequency, {$\jmath$ denotes the imaginary unit},
$\epsilon_0$ and $\mu_0$ are the absolute permittivity and permeability of vacuum, respectively. Note the presence of the artificial non-homogeneous and non-isotropic material,  {$\bar{\bm\epsilon}_r$, $\bar{\bm\mu}_r$ being the relative permittivity and permeability in \eqref{permit}-\eqref{permea} originated by the coordinate transformation \eqref{rrp}-\eqref{eq:transformation}}.

The vector wave equation for the electric field ${\bf{E}}$ derived from \eqref{eq:maxwell1} and \eqref{eq:maxwell2} is 
\begin{equation}\label{eqwave}
\nabla\times\bar{\bm\mu}_r^{-1}\nabla\times{\bf E}-k_0^2\bar{\bm\epsilon}_r{\bf E} = 0,
\end{equation}
where $k_0=\omega\sqrt{\mu_0\epsilon_0}$ is the free space wavenumber. As discussed previously, a perfect conductor boundary condition is imposed in an essential way on the lateral surface, $S_L$, of the taper, while the magnetic field is assigned as a natural condition on the device ports, $S_1$ and $S_2$.
With a view to the HiMod discretization, we consider the weak form of problem \eqref{eqwave}, i.e.,  we look for a solution ${\bf E} \in {H_{S_L}({\rm curl},\Omega)}$, such that
\begin{equation}
\begin{split}
&\int_{\Omega}(\nabla\times{\bf w})\cdot\bar{\bm\mu}_r^{-1}\nabla\times{\bf E}\;d\Omega-k_0^2\int_\Omega{\bf w}\cdot\bar{\bm\epsilon}_r{\bf E}\;d\Omega\\
&-{\int_{S_1\cup S_2}}{\bf w}\cdot\bar{\bm\mu}_r^{-1}[(\nabla\times{\bf E})\times\hat{\bf n}]\;dS = 0,
\end{split} \label{eqweak}
\end{equation}
with $\hat{\bf n}$ the unit outward normal vector to $\partial\Omega$, for any test function, ${\bf w}$, in ${H_{S_L}({\rm curl},\Omega)}$, the space of the functions that are square integrable in $\Omega$ together with the associated curl, $\nabla\times{\bf w}$~\cite{Nedelec}{, satisfying $\nabla\times{\bf E}\times\hat{\bf n}={\bf 0}$ on $S_L$.}

Now, to simplify the notation, 
we no longer distinguish between the transverse modal fields of type TE or TM {in the HiMod expansion \eqref{eq:et0}}, and we order the modes by using a single index, $n$, that runs from $1$ to $N_{\rm M} = N_{\rm TE}+N_{\rm TM}$. Hence, the electric field expansions in \eqref{eq:et0}-\eqref{eq:ez0}, {taking into account the FE representations in \eqref{eqtetau}-\eqref{eqtmz}}, can be rewritten as
\begin{align}\label{eqEt}
&{\bf{E}}_t(x,y,z) = \sum_{l=1}^{N_{l,t}}\sum_{n=1}^{N_{\rm M}}c_{n,l}{\bf{e}}_{t,n}(x,y)\phi_l(z) \\ \label{eqEz}
&{\bf{E}}_z(x,y,z) = \sum_{l=1}^{N_{l,z}}\sum_{n=1}^{N_{\rm TM}}d_{n,l}{\bf{e}}_{z,n}(x,y)\psi_l(z).
\end{align}
The following normalization is then applied to the transverse and to the longitudinal fields, so that
\begin{equation}\label{eqnorm0}
{\int_{S_P}} {\bf e}_{t,n}\cdot{\bf e}_{t,n}\;dS = 1,\hspace{1em}
{\int_{S_P}} {\bf e}_{z,n}\cdot{\bf e}_{z,n}\;dS = 1,
\end{equation}
{where $S_P$ denotes the port surface after the coordinate transformation, with $P\in\{1, 2\}$.} 

{By arranging the separate modal and FE basis functions in \eqref{eqEt}, we can define the generic $i$-th global basis function involved in the expansion of the transverse component of the electric field 
	as:}
\begin{equation}\label{f_i}
{\bf f}_{i}(x,y,z) = {\bf e}_{t,n_i}(x,y)\phi_{l_i}(z)
\end{equation}
{where the sorting is made so that
	$(n_1,l_1)=(1,1), (n_2,l_2)=(2,1), \ldots,
	(n_{N_{\rm M}},l_{N_{\rm M}})=(N_{\rm M},1)$,
	$(n_{N_{\rm M}+1},l_{N_{\rm M}+1})=(1,2), \ldots,
	(n_{N_{\rm M}\cdot N_{l,t}},l_{N_{\rm M}\cdot N_{l,t}})=(N_{\rm M},N_{l,t})$.
	
	Analogously, from \eqref{eqEz}, the generic $i$-th {global} basis function} for the longitudinal component of field ${\bf E}$ 
takes the form:
\begin{equation}\label{g_i}
{\bf g}_{i}(x,y,z) = {\bf e}_{z,n_i}(x,y)\psi_{l_i}(z),
\end{equation}
{where, now, the index sorting is such that
	$$
    (n_1,l_1)=(1,1), (n_2,l_2)=(2,1), \ldots,
	 (n_{N_{\rm TM}},l_{N_{\rm TM}})=(N_{\rm TM},1),
	 $$
	$$(n_{N_{\rm TM}+1},l_{N_{\rm TM}+1})=(1,2), \ldots,
	(n_{N_{\rm TM}\cdot N_{l,z}},l_{N_{\rm TM}\cdot N_{l,t}})=(N_{\rm TM},N_{l,z}).
	$$
We remark that the index notation $n_i$ and $l_i$ 
{in \eqref{f_i}-\eqref{g_i}} is meant to highlight that the generic $i$-th {global} basis function is associated with a modal index, $n_i$, {(which, in turn, depends on indices $p_{n_i}$ and $q_{n_i}$ according to \eqref{pnindices}) related to the approximation of the transverse dynamics, and with a FE index, $l_i$, related to} the discretization along $z$.

{By exploiting the separation of variables characterizing the definition of functions ${\bf f}_{i}$ and ${\bf g}_{i}$, we have
}

\begin{align}\label{eqcomp1}
\nabla\times{\bf f}_{i} &= \phi_{l_i}\nabla_t\times{\bf e}_{t,n_i}+ \frac{d\phi_{l_i}}{dz}\hat{\bf z}\times{\bf e}_{t,n_i}, \\ \label{eqcomp2}
\nabla\times{\bf g}_{i} &= \psi_{l_i}\nabla_t\times{\bf e}_{z,n_i},
\end{align}
with $\nabla_t$ the transverse component of operator $\nabla$ in the $xy$-plane, and $\hat{\bf z}$ the unit vector along the $z$-direction.\\
For the compatibility of the representations of the transverse component of the curl operator, the degree of polynomials $\psi_{l_i}$ has to be one less than the degree of polynomials $\phi_{l_i}$. In fact, the second term in the right-hand side of \eqref{eqcomp1} and the term in the right-hand side of \eqref{eqcomp2} both define the transverse components of the curl and different degrees in the polynomial representations lead to spurious solutions~\cite{GentiliEplane,Cendes}.

{Now, we assume to deal with a generic excitation ${\bf h}$, which represents a magnetic field impressed at the ports $S_1$ and $S_2$. 
	The HiMod discretization in \eqref{eqEt} and \eqref{eqEz}, combined with definitions \eqref{f_i}-\eqref{g_i} and relations \eqref{eqcomp1}-\eqref{eqcomp2},} convert problem \eqref{eqweak} into the following matrix formulation
\begin{equation}\label{eqmatAB}
\left(
\begin{bmatrix}
{\bm A}_{tt} & {\bm A}_{tz}\\
{\bm A}_{zt} & {\bm A}_{zz}\\
\end{bmatrix}
- k_0^2
\begin{bmatrix}
{\bm B}_{tt} &  {\bm B}_{tz}\\
{\bm B}_{zt} & {\bm B}_{zz}\\
\end{bmatrix}
\right) 
\begin{bmatrix}
{\bm c} \\ {\bm d}
\end{bmatrix}
= -\jmath\omega\mu_0
\begin{bmatrix}
{\bm C}_t \\ {\bf 0}
\end{bmatrix},
\end{equation}
where vectors
$$
\begin{array}{rcl}
{\bm c} &=& \big[c_{n_1,l_1}, \ldots, c_{n_{N_{\rm M}\cdot N_{l,t}}, \,l_{N_{\rm M}\cdot N_{l,t}}}\big]^T,\\[2mm]
{\bm d} &=& \big[d_{n_1,l_1}, \ldots, d_{n_{N_{\rm TM}\cdot N_{l,z}},\, l_{N_{\rm TM}\cdot N_{l,t}}}\big]^T
\end{array}
$$
collect the (unknown) HiMod coefficients in \eqref{eqEt}-\eqref{eqEz}, while 
the generic matrix entries are given by
\begin{align}
&{\bm A}_{tt}(i,j) =\int_{\Omega}{\nabla\times{\bf f}_i\cdot\bar{\bm\mu}_r^{-1}\nabla\times{\bf f}_j}\,d\Omega, \label{att} \\
&{\bm A}_{tz}(i,j) =\int_{\Omega}{\nabla\times{\bf f}_i\cdot\bar{\bm\mu}_r^{-1}\nabla\times{\bf g}_j}\,d\Omega,  \\
&{\bm A}_{zt}(i,j) =\int_{\Omega}{\nabla\times{\bf g}_i\cdot\bar{\bm\mu}_r^{-1}\nabla\times{\bf f}_j}\,d\Omega, \\
&{\bm A}_{zz}(i,j) =\int_{\Omega}{\nabla\times{\bf g}_i\cdot\bar{\bm\mu}_r^{-1}\nabla\times{\bf g}_j}\,d\Omega, 
\end{align}
and
\begin{align}
&{\bm B}_{tt}(i,j) =\int_{\Omega}{{\bf f}_i\cdot\bar{\bm\epsilon}_r{\bf f}_j}\,d\Omega, \\
&{\bm B}_{tz}(i,j) =\int_{\Omega}{{\bf f}_i\cdot\bar{\bm\epsilon}_r{\bf g}_j}\,d\Omega, \\
&{\bm B}_{zt}(i,j) =\int_{\Omega}{{\bf g}_i\cdot\bar{\bm\epsilon}_r{\bf f}_j}\,d\Omega, \\
&{\bm B}_{zz}(i,j) =\int_{\Omega}{{\bf g}_i\cdot\bar{\bm\epsilon}_r{\bf g}_j}\,d\Omega.
\end{align}
Finally, vector ${\bm C}_t$ on the right-hand side in \eqref{eqmatAB} is defined by
\begin{equation}\label{eqCmat}
{{\bm C}_t(i)} = {\int_{S_1\cup S_2}}{\bf f}_i\cdot({\bf h}\times\hat{\bf n})\,dS,
\end{equation}
where we set $\bar{\bm\mu}_r^{-1}\nabla\times{\bf E} = -\jmath\omega\mu_0{\bf h}$. 
Gathering matrices and vectors in \eqref{eqmatAB} as
\begin{equation}
{\bm A} = 
\begin{bmatrix}
{\bm A}_{tt} & {\bm A}_{tz}\\
{\bm A}_{zt} & {\bm A}_{zz}\\
\end{bmatrix},\quad
{\bm B} = 
\begin{bmatrix}
{\bm B}_{tt} & {\bm B}_{tz}\\
{\bm B}_{zt} & {\bm B}_{zz}\\
\end{bmatrix},
\end{equation}
\begin{equation}\label{eqC}
{\bm C} = 
\begin{bmatrix}
{\bm C}_t \\
{\bf 0}
\end{bmatrix},
\hspace{2em}
{\bm v} = 
\begin{bmatrix}
{\bm c}\\
{\bm d} \\
\end{bmatrix},
\end{equation}
yields {the HiMod solution 
	\begin{equation}\label{vsolution}
	{\bm{v}} = -\jmath\omega\mu_0({\bm A}-k_0^2{\bm B})^{-1}{\bm C}
	\end{equation}
	to equation \eqref{eqweak}.}

\subsection{Generalized scattering matrix computation} 
The magnetic field ${\bf h}$ is characterized by a different expression at port $1$ and at port $2$. Indeed, by construction, port $2$ is transformed through the mapping in \eqref{eq:transformation}, while port $1$ is not. 
Now, we compute the Generalized Scattering Matrix (GSM) starting from the normalized Generalized Impedance Matrix (GIM).
We denote by ${\bf e}^{(P)}_{n}$, ${\bf h}^{(P)}_{n}$ the electric and the magnetic modal functions of the $z$-propagating mode of index $n$ at port $P\in\left\{1,2\right\}$,  normalized as
\begin{equation}\label{eqnorm}
\int_{S_P}{\bf e}^{(P)}_{n}\times{\bf h}^{(P)}_{n}\cdot\hat{\bf z}\;dS = \pm 1,
\end{equation}
where signs '$+$' and '$-$' refer to port $1$ and $2$, respectively (oppositely oriented with respect to the $z$-axis). Relation \eqref{eqnorm} is the standard normalization employed in power waves, commonly used to define the scattering matrix.

We now specify the actual expression of modal functions ${\bf e}^{(P)}_{n}$ and ${\bf h}^{(P)}_{n}$ at the two ports, separately. 

\subsubsection{Port $1$}
For this port, modes are the same functions used to expand the electric field. Therefore, by using \eqref{eqnorm0} and enforcing \eqref{eqnorm}, we have
\begin{align}
&{\bf e}_n^{(1)} = \frac{{\bf e}_{t,n}}{\sqrt{Y_{n}^{(1)}}} \\
&{\bf h}_n^{(1)} = {\hat{\bf z}\times{\bf e}_{t,n}}{\sqrt{Y_{n}^{(1)}}},
\end{align}
being 
\begin{equation}\label{key}
Y_{n}^{(1)} = \left\{
\begin{matrix}
\cfrac{\gamma_{n}^{(1)}}{\jmath\omega{\mu'}}\hspace{1em}\text{\rm for TE modes} \\[3mm]
\cfrac{\jmath\omega{\epsilon'}}{\gamma_{n}^{(1)}}\hspace{1em}\text{\rm for TM modes}		
\end{matrix}
\right. 
\end{equation}
the modal wave admittance at port $1$, where $\gamma_{n}^{(1)} = \sqrt{[{k_{n}^{(1)}}]^2-k^2}$ is the modal propagation constant, with $k_{n}^{(1)}$ the mode eigenvalue (known from the standard modal decomposition in rectangular waveguides) and $k=\omega\sqrt{\mu'\epsilon'}$ the wavenumber, being $\varepsilon'$, $\mu'$ the absolute permittivity and permeability of the waveguide at port 1. 

\vspace*{0.2cm}

\subsubsection{Port $2$}
This port is subject to the geometric transformation in \eqref{eq:transformation} so that sides of length $a_L$, $b_L$ are transformed into sides with length $a_0$, $b_0$. The modal functions are subject to transformation 
\eqref{eq:transformation} for $z'=L$, namely,
\begin{equation}\label{eqtrasf1}
x = x'\frac{a_0}{a_L}, \quad
y = y'\frac{b_0}{b_L}, \quad
z = z',
\end{equation}
so that the associated Jacobian matrix, $\bar{\bf J}_L$, affecting the output port is 
\begin{equation}
\bar{\textbf{J}}_L = \begin{bmatrix}
\frac{a_0}{a_L} & 0 & 0 \\
0 & \frac{b_0}{b_L} & 0 \\
0 & 0 & 1 		
\end{bmatrix}.
\end{equation}
It can be shown that the propagation constants, $\gamma_{n}^{(2)}$, of the modes at port $2$ are not modified by the transformation, while the modal functions are transformed, i.e., we have
\begin{align}
&{\bf e}_n^{(2)} = A_n\;\bar{\bf j}^{-1}_L{\bf e}_{t,n},\label{eqhtr0} \\[2mm]
\label{eqhtr1}
&{\bf h}_n^{(2)} = -A_nY_{n}^{(2)}\;\bar{\bf j}^{-1}_L(\hat{\bf z}\times{\bf e}_{t,n}),
\end{align}
where
\begin{equation}
\bar{\textbf{j}}_L = \begin{bmatrix}
\frac{a_0}{a_L} & 0 \\
0 & \frac{b_0}{b_L} 
\end{bmatrix},
\end{equation}
\begin{equation}
Y_{n}^{(2)} = \left\{
\begin{matrix}
\cfrac{\gamma_{n}^{(2)}}{\jmath\omega\mu'}\hspace{1em}\text{\rm for TE modes} \\[3mm]
\cfrac{\jmath\omega\epsilon'}{\gamma_{n}^{(2)}}\hspace{1em}\text{\rm for TM modes}		
\end{matrix}
\right. 
\end{equation}
is the modal wave admittance at port $2$, with $\gamma_{n}^{(2)} = \sqrt{[{k_{n}^{(2)}}]^2-k^2}$ {the modal propagation constant, $k_{n}^{(2)}$ the mode eigenvalue and $k$ defined as in \eqref{key}}, and $A_n$ is a normalization constant {to be properly set. In particular,} imposing \eqref{eqnorm} leads to
\begin{equation*}
A_n = \sqrt{\frac{a_0b_0}{Y_{n}^{(2)}a_Lb_L}},
\end{equation*}
so that modal functions in \eqref{eqhtr0}-\eqref{eqhtr1} become
\begin{align}
&{\bf e}_n^{(2)} = \sqrt{\frac{a_0b_0}{Y_{n}^{(2)}a_Lb_L}}\;\bar{\bf j}^{-1}_L{\bf e}_{t,n}, \\
&{\bf h}_n^{(2)} = -\sqrt{\frac{Y_{n}^{(2)}a_0b_0}{a_Lb_L}}\;\bar{\bf j}^{-1}_L(\hat{\bf z}\times{\bf e}_{t,n}).
\end{align}

\vspace*{0.2cm}

We now detail the computation of matrix GSM. The idea is to modify the vector ${\bm C}$ in \eqref{eqC} into a matrix in order to include all the modal functions at the two ports, namely
\begin{equation}
{\bm C} = 
\begin{bmatrix}
{\bm C}_t^{(1)} & {\bm C}_t^{(2)} \\[1mm]
{\bf 0} & {\bf 0}
\end{bmatrix},
\end{equation}
where we use again the superscript ${(P)}$ to refer to port $P\in\{1,2\}$, and
\begin{equation}
{\bm C}_t^{(P)}(i, n) =\int_{S_P}{\bf f}_i\cdot[{\bf h}_n^{(P)}\times\hat{\bf n}]\;dS.
\end{equation}

The expression for GIM, ${\bm Z}$, is finally obtained from \eqref{vsolution} with a multi-column forcing term ${\bm C}$, which yields
\begin{equation}\label{ZZ}
{\bm Z} = -j\omega\mu_0{\bm C}^T({\bm A}-k_0^2{\bm B})^{-1}{\bm C}.
\end{equation}
It is well known \cite{collin2007foundations} that matrix ${\bm Z}$ defines the relationship between the vector, ${\bm p}^+$, of incident power waves and vector, ${\bm p}^{-}$, of the reflected power waves, for each mode at each port, being
\begin{equation}
{\bm p}^++{\bm p}^- = {\bm Z}\;({\bm p}^+-{\bm p}^-).
\end{equation}
Thus, setting ${\bm p}^-={\bm S}{\bm p}^+$, one can compute
\begin{equation}\label{eqS}
{\bm S} = ({\bm Z}+{\bm I})^{-1}({\bm Z}-{\bm I}),
\end{equation}
with ${\bm I}$ the identity matrix of size $2N_{\rm M}\times 2N_{\rm M}$, which provides the standard expression for the scattering matrix. 

Note that we have defined all terms needed to build the final GSM, ${\bm S}$,  by using the transverse modal functions ${\bf e}_{t,n}$ only.

\section{Results}\label{sec4}
In this section we assess the performance of the HiMod discretization when employed to compute
the GSM for different instances of smooth transitions in a rectangular waveguide. We investigate 
the accuracy and the computational efficiency of the proposed method, when compared with a standard full 3D FEM. In particular, we adopt the results obtained with the commercial software CST Microwave Studio (version $18$) and HFSS (Ansys Electronic Desktop, Release 2021 R2) as reference values, by pushing the accuracy to get very accurate data. 
In particular, we adopt a frequency-domain iterative solver with tolerance equal to $1$e$-4$ for CST and $1$e$-5$ for HFSS.
\begin{figure}[!t]
	\centering
	\includegraphics[width=\columnwidth]{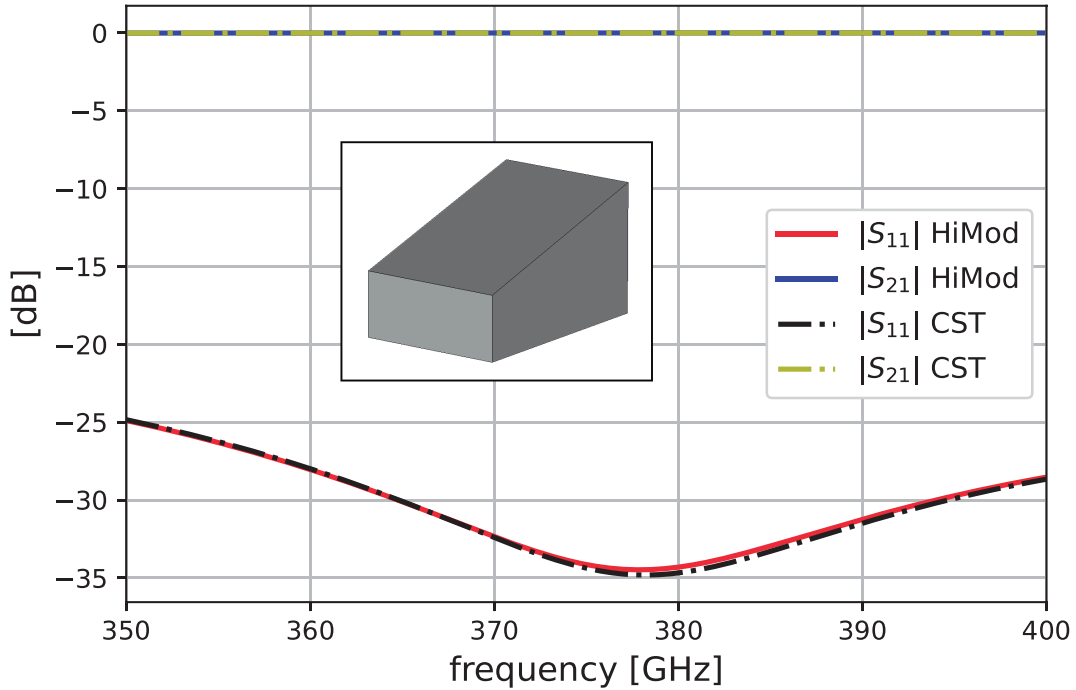}
	\caption{
		Magnitude of the scattering parameters of the TE\textsubscript{10} mode for a taper with linear variation in height as in~\cite{Garner}. Comparison between results with commercial software CST and the HiMod reduced solution.}
	\label{fig:easy}
\end{figure}

As a first example, we analyze {a linear taper which connects two waveguides characterized by a rectangular cross-section, varying in height only (see Fig.~\ref{fig:easy}). In particular, following~\cite{Garner}, we choose
	$a_0=0.570$ mm, $b_0=0.570$ mm, $a_L=0.285$ mm, $b_L=0.570$ mm, $L=1.1$ mm $= 1.36\lambda_0$,
	being $\lambda_0=c/f_0$ with $c$ the light speed and $f_0=372$ GHz, according to the rule
	$$
	a(z) = a_0+(a_L-a_0)\frac{z}{L}, \quad b(z) = b_0.
	$$
	In~\cite{Garner}, this simple taper is used as a component in a higher-order mode input coupler at $372$ GHz.\\
	For this very simple filter, we 
	discretize the direction $z$ by resorting to quadratic finite elements based on
	$11$ DOFs only, while $4$ modes (${\rm TE}_{10}$, ${\rm TE}_{01}$, ${\rm TE}_{11}$ and ${\rm TM}_{11}$) are adopted to model the transverse field, amounting to just $50$ unknowns for each frequency point employed to solve \eqref{ZZ}. Components \eqref{att}-\eqref{eqCmat} are computed by exploiting the iterated integral rule thanks to the Cartesian 
	geometry of the transformed taper. In particular, we use a
	1D Gauss-Legendre rule along $x$, $y$ and $z$ with a relative tolerance equal to $1.5$e$-05$.
	
	In Fig.~\ref{fig:easy}, we compare the results provided by CST software and the HiMod procedure in terms of parameters $S_{11}$ (the reflection coefficient at port $1$) and $S_{21}$ (the transmission coefficient at port $2$) for the dominant TE$_{10}$ mode in the input/output waveguides. Other non-propagating modes in the bandwidth of analysis are present at the ports, but the plots refer to the dominant TE$_{10}$ mode only to ease the comparison with some reference data. \\
	The curves in the figure exhibit 
	a very good matching between the CST and the HiMod solutions. 
	
	As a more general configuration, we choose a waveguide with dimensions $a_0=22.86$ mm, $b_0=11.43$ mm, that is connected to a bigger one, with $a_L=28.448$ mm, $b_L=14.224$ mm, by means of a linear taper of length $L=20$ mm $=0.67\lambda_0$ at $f_0=10$ GHz (see Fig.~\ref{fig:linear}). This second setting introduces a change of size along both the $x$- and $y$-direction. The profile of the device is analytically defined by the transformation
	\begin{align}
	&a(z) = a_0+(a_L-a_0)\frac{z}{L}, \\[2mm]
	&b(z) = b_0+(b_L-b_0)\frac{z}{L}.
	\end{align}
	The finite element and the modal discretizations have been selected to guarantee a sufficient accuracy to the HiMod approximation in order to match the CST solution. In particular, we discretize the $z$-direction with $15$ equally-spaced points, for a total of $15$ linear and $29$ quadratic DOFs.
	Concerning the modal expansion, we employ the same $4$ modes as in the previous case.\\
	In  Fig. \ref{fig:linear}, we compare CST and HiMod solutions for the dominant mode TE$_{10}$. 
	The results are essentially identical for the entire bandwidth.
	\begin{figure}[!t]
		\centering
		\includegraphics[width=\columnwidth]{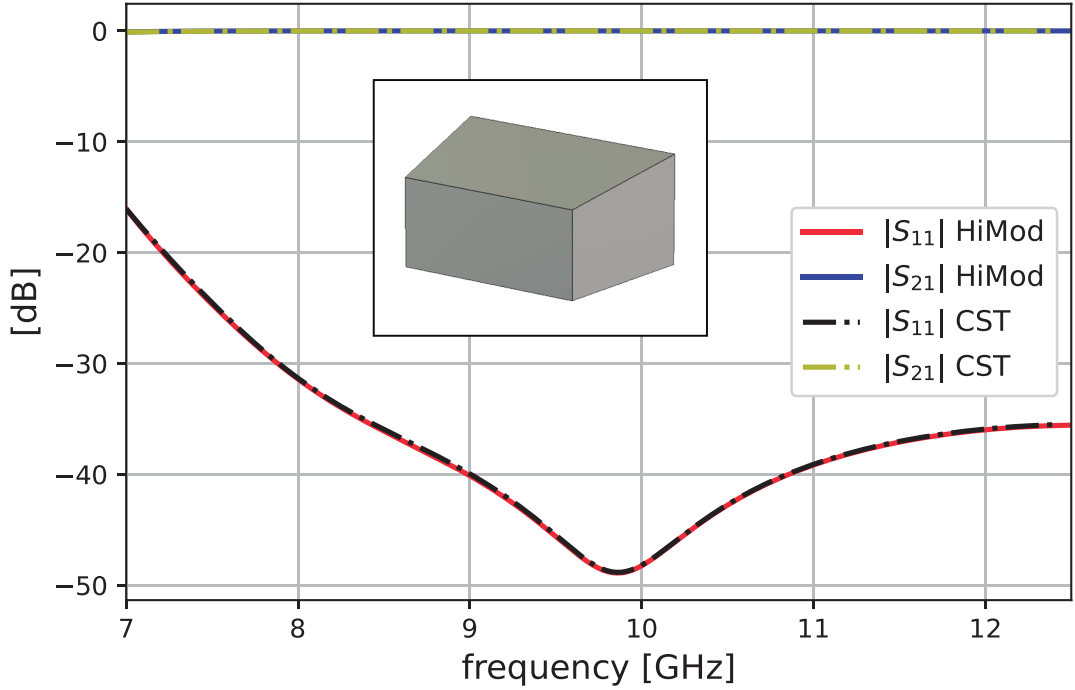}
		\caption{Magnitude of the scattering parameters of the TE\textsubscript{10} mode for a taper with linear variation in  height and width. Comparison between results with commercial software CST and the HiMod reduced solution.}
		\label{fig:linear}
	\end{figure}
	
	The third example refers to the sinusoidal taper
	shown in Fig.~\ref{fig:sin}, where the waveguide dimensions change only along the $x$-direction. The taper is now described by relations
	$$a(z) = a_0+ (a_L-a_0)\sin\left(\frac{\pi z}{2L}\right), \quad b(z) = b_0.
	$$
	where the waveguide dimensions are $a_0=15.79$ mm, $b_0 = 7.889$ mm, and $a_L=22.86$ mm, $b_L = 7.889$ mm, with a length $L=40$ mm $= 1.67\lambda_0$ at $f_0=12.5$ GHz.\\ 
	The HiMod reduction is based on the same discretization along $z$ as for the previous case study, whereas the four TE\textsubscript{m0} ($m=1,...,4$) modes are employed in the transverse section, without any TM mode. This leads to solve an algebraic system of size equal to $116$.\\
	In Fig. \ref{fig:sin}, the response for the first propagating mode TE$_{10}$ covers higher
	frequencies when compared with the ones in Fig.~\ref{fig:linear}. The results associated with the HiMod reduced model are very accurate and in excellent agreement with the output provided by CST, except for the very small mismatch for the highest values of the frequency range, corresponding anyway to an extremely small value of reflection coefficient.
	\begin{figure}[!t]
		\centering
		\includegraphics[width=\columnwidth]{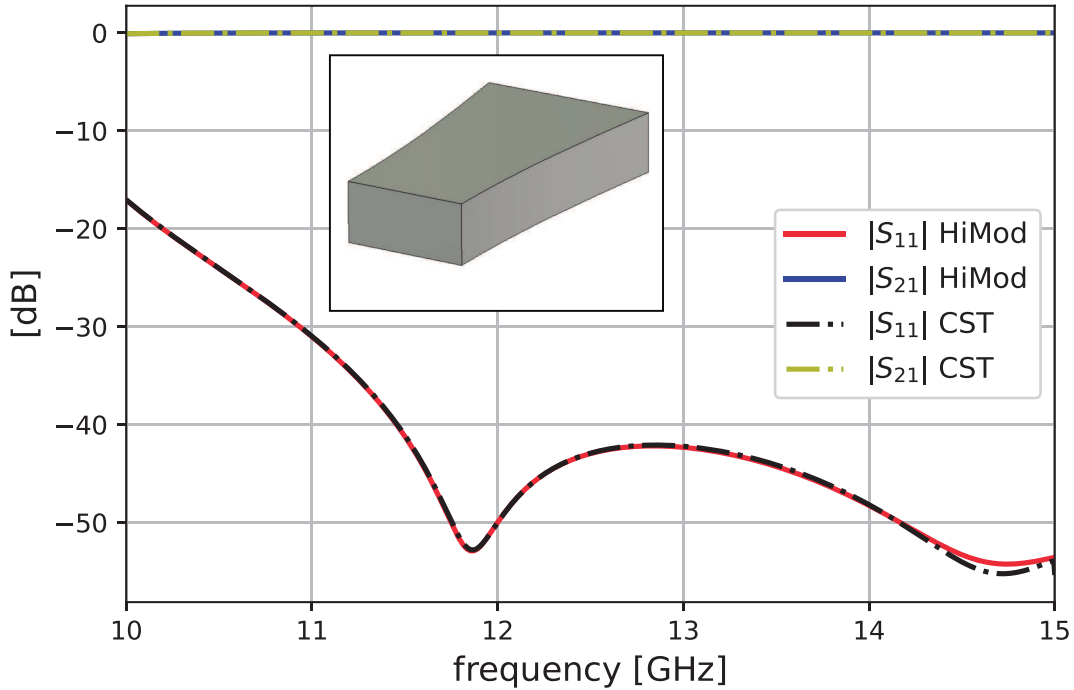}
		\caption{Magnitude of the scattering parameters of the TE\textsubscript{10} mode for a taper with sinusoidal variation 
			along the $x$-direction. Comparison between results with
			commercial software CST and the HiMod reduced solution.}
		\label{fig:sin}
	\end{figure}
	
	As a practical example taken from the literature, we simulate the taper shown in Fig.~\ref{fig:hecken}  by setting $a_0=22.86$ mm, $b_0 = 10.16$ mm, $a_L=22.86$ mm, $b_L = 5.08$ mm and length $L=47.08$ mm $=1.41\lambda_0$, being $f_0=10$ GHz.
	The profile synthesis was done by the procedure described in~\cite{Percaz} for a Hecken taper \cite{Hecken}, which exhibits a quasi-optimal response for the given set of input parameters (i.e., a return loss lower than $40$ dB for the WR90 frequency range, between $8.2$ GHz and $12.4$ GHz). The $b(z)$ dependence does not have an explicit form but can be expressed as
	\begin{equation}\label{Kb}
	b(z) = b_0\exp\left({-2\int_{0}^{z}\!K_b(r)dr}\right),
	\end{equation}
	where $K_b(r)$ is a function of the return loss, the frequency range, and of the dimension of both the ports, while $a(z)=a_0$. The integral involved in definition \eqref{Kb} is, in general, computed numerically. Interested readers are referred to~\cite{Percaz} for further details, as the process cannot be easily summarized.\\ 
	We employ the HiMod approximation by preserving the discretization along $z$ as in the two previous test cases, and exploiting modes TE\textsubscript{10}, TE\textsubscript{12}, TE\textsubscript{14} and TM\textsubscript{12}, TM\textsubscript{14} to model the transverse fields, for a total number of $175$ unknowns.\\
	The results shown in Fig.~\ref{fig:hecken} highlight the 
	good accuracy characterizing the reduced solution when compared with the CST output.
	The slight difference in the minimum of the second resonance is likely due to the frequency step employed to compute the GSM, and does not represent a crucial issue. 
	In such a context, it is not trivial to establish which method is better.  
	Indeed, matching results of the order of $-60$ dB are very satisfactory, showing that the way we are modelling tapers with HiMod does not incur in any appreciable simplification of the geometry.
	\begin{figure}[!t]
		\centering
		\includegraphics[width=\columnwidth]{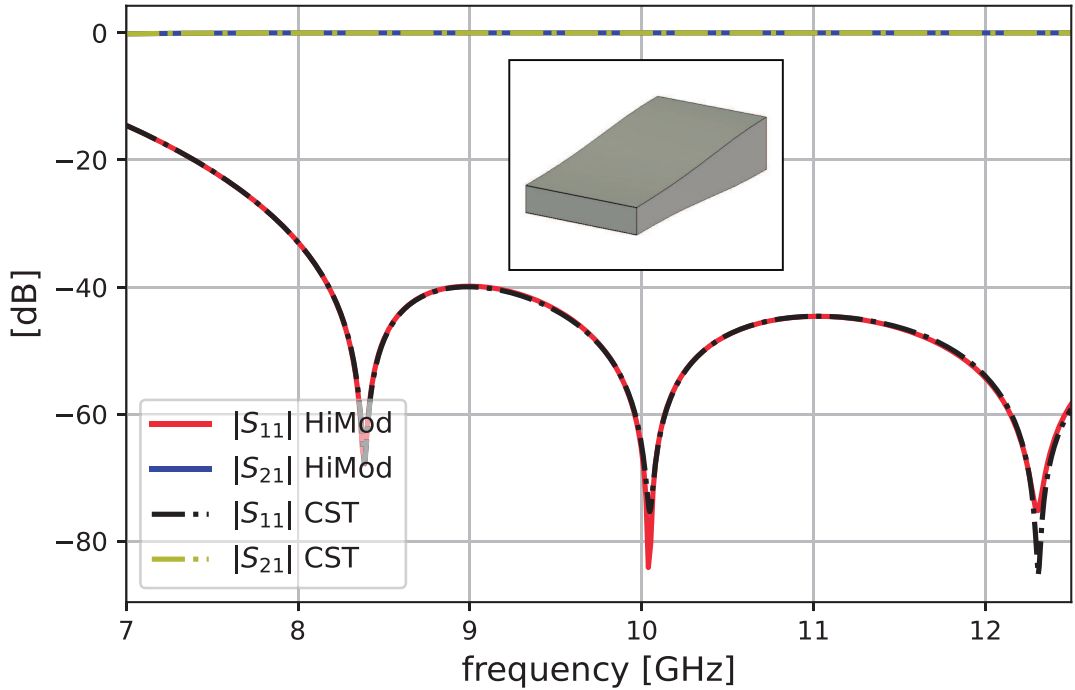}
		\caption{Magnitude of the scattering parameters of the TE\textsubscript{10} mode in one of the Hecken tapers in \cite{Percaz}. Comparison between results with commercial software CST and the HiMod reduced solution.}
		\label{fig:hecken}
	\end{figure}
	
	The numerical assessment carried out in the previous examples can be conceived as a benchmark case study due to the low number of DOFs involved by both the 3D FEM and the HiMod discretizations.
	In order to challenge the reduced order model in terms of accuracy and scalability, 
	we simulate an E-plane low-pass filter built on a WR-75 waveguide, with input dimensions $a_0=a_L=19.05$ mm, $b_0=b_L= 9.525$ mm, $L=218.025$ mm $= 5.45\lambda_0$ at $f_0=12$ GHz, composed by sinusoidal units of length $3.825$ mm and a minimum gap of $6$ mm, as proposed in~\cite{Arregui2010, Arregui2013}.
	
	The filter profile is shown in Fig. \ref{fig:filter} and is made up with $3$ sections. The first one consists of 
	$8$ units modulated by a Hanning window of height (peak-to-peak) $h_{min}=1.6$ mm. The second section comprises $45$ units of increasing height, from $h_{min}$ to $h_{max}=8.6$ mm, the last four being set to $h_{max}$, and with a minimum gap equal to $g_{min}= 6$ mm. The final section is constituted by $4$ units with a maximum height given by $h_{max}$, and are modulated by a Hanning window (for more details, we refer the reader to~\cite{Arregui2010}).\\
	We observe that the smooth variation of the boundaries of the device makes it especially suitable for high-power applications. However, from a computational viewpoint, 
	the high frequency variation of the sinusoidal profile 
	turns out to be very challenging for 3D methods. This justifies the adoption of the HiMod reduced model as well as of the dedicated E-plane 2D-FEM formulation for rectangular waveguides proposed in~\cite{GentiliEplane}.
	In particular, the 2D formulation employs quadratic finite elements to approximate the $x$-component of the magnetic  field, with about $400,000$ DOFs. The HiMod reduction 
	discretizes the propagating direction with $451$ equally space points, the transverse section with $7$ modes (TE\textsubscript{10}, TE\textsubscript{12}, TM\textsubscript{12}, TE\textsubscript{14}, TM\textsubscript{14}, TE\textsubscript{16} and TM\textsubscript{16}), with a total number of DOFs equal to {$7660$}, after subdividing the frequency range $10-15$ GHz with $201$ samples (about $0.15$s per sample). 
	
	In Fig.~\ref{fig:sparam_filter}, the filter response has been modeled by comparing HiMod approximation with that obtained by the commercial software HFSS (which simplifies the set-up of this configuration with respect to CST) and by the 2D-FEM model. 
	We remark that
	the three numerical procedures accurately capture the response in the passband. However, only HiMod and the 2D-FEM formulation (which is a less general procedure than HiMod in the context of this work) succeed in computing the extremely low magnitude of paramater $S_{21}$ in the stopband, virtually achieving full rejection. Moreover, among the three approaches,
	HiMod requires the least amount of unknowns and is the fastest one. 
	From a more quantitative viewpoint, the HFSS simulation takes over three hours on a standard PC, while the 2D-FEM and the HiMod simulation require three minutes and less than one minute, respectively, without any code parallelization, only using sparse matrices and open-source solvers.
	\begin{figure}[!t]
		\centering
		\includegraphics[width=\columnwidth]{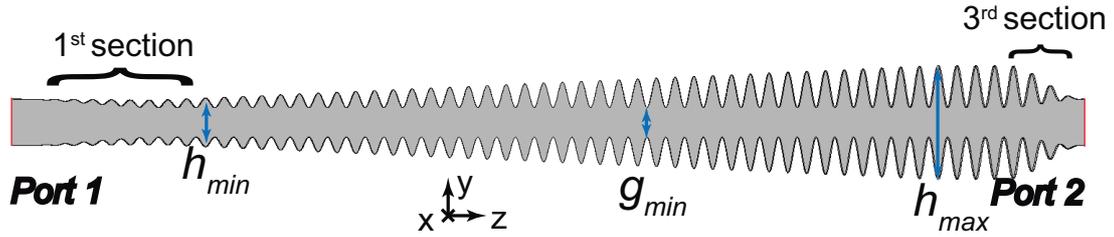}
		\caption{Rectangular-waveguide E-plane low-pass filter ($y$-$z$ view) described in \cite{Arregui2010}.}
		\label{fig:filter}
	\end{figure}
	
	\begin{figure}[!t]
		\centering
		\includegraphics[width=\columnwidth]{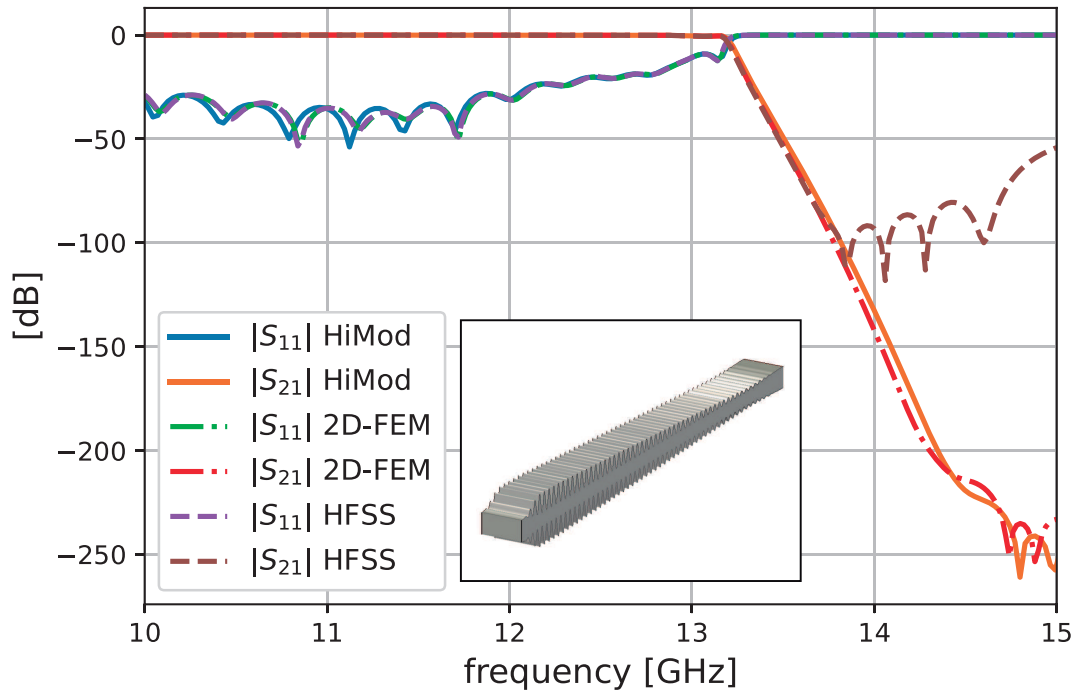}
		\caption{Response of the of the TE\textsubscript{10} mode in filter in Fig. \ref{fig:filter} for HFSS, a dedicated E-plane 2D-FEM and HiMod.}
		\label{fig:sparam_filter}
	\end{figure}
	
	\section{Conclusions}\label{sec5}
	A new method based on a Hierarchical Model (HiMod) reduction has been presented and applied to the analysis of rectangular waveguide tapers and transitions with a smooth profile. The method is based on a preliminary coordinate transformation that converts the device into a parallelepiped filled by an inhomogeneous anisotropic material. In the transformed geometry, a modal decomposition in the transverse plane combined  with  1D finite elements along the longitudinal direction is then adopted to model the electromagnetic wave propagation. 
	\\
	{We have verified that the actual 
		advantage of the HiMod reduction resides in the very small size of the associated algebraic system, without waving the accuracy. 
		In the most challenging considered configuration (see Fig.~\ref{fig:filter}), $7$ modes are sufficient to guarantee an accurate reduced solution which outperforms the output provided by HFSS software. Moreover, extremely accurate results have been obtained with an impressive speed-up factor (about $200$) when compared to standard 3D FEM, although no optimization has been applied to the code of the reduced model.}
	
	As a possible future development, we are interested in varying the number of modes along the $z$-direction according to a pointwise HiMod formulation (see~\cite{PerottoZilio13,PerottoZilio15}).
	
	\section*{Acknowledgement}
	This work was supported in part by the Spanish Government under grant PID2020-116968RB-C32 (DEWICOM) and TED2021-130650B-C21 (ANT4CLIM) funded by MCIN/ AEI/ 10.13039/501100011033 (Agencia Estatal de Investigación) and by UE (European Union) ``NextGenerationEU"/PRTR. S. Perotto acknowledges the European Union's Horizon 2020 research and innovation programme under the Marie Skłodowska-Curie Actions, grant agreement 872442 ({\it ARIA, Accurate Roms for Industrial Applications}). S. Micheletti and S. Perotto thank the PRIN research grant n.20204LN5N5 {\it Advanced Polyhedral Discretisations of Heterogeneous PDEs for Multiphysics Problems} and the INdAM-GNCS 2022 Project {\it Metodi di riduzione computazionale per le scienze applicate: focus su sistemi complessi}.

	\appendix
	
	\section{Modal functions in the rectangular waveguide}
		
		Here, we define the expressions for functions ${\bf e}_{t,n}^{\rm TX} = \hat{\bf x}e_{x,n}^{\rm TX}+\hat{\bf y}e_{y,n}^{\rm TX}$, with ${\rm TX}={\rm TE}$, ${\rm TM}$ and ${\bf e}_{z,n}^{\rm TM} = \hat{\bf z}e_{z,n}^{\rm TM}$ in \eqref{eq:et0}-\eqref{eq:ez0}, which are the transverse and longitudinal modal functions, respectively associated with a generic mode $n$.\\ 
		Modes in the rectangular waveguide are defined by two indexes $p = p(n)$ and $q = q(n)$, associated with the $x$ and $y$ coordinates, respectively (see \eqref{pnindices}) and spanning the ranges $p = 0, 1\ldots$ and $q = 0, 1, \ldots$, except for the choice $p=q=0$. \\As discussed in Section \ref{sec2}, we consider a waveguide whose cross section is centered in the $x-y$ plane, of size $a_0$ and $b_0$ along $x$ and $y$, respectively. For $n=1, \ldots, N_{\rm M}$, with $N_{\rm M} = N_{{\rm TE}}+N_{{\rm TM}}$, if the $n$-th mode is TE, it holds
		\begin{align}\label{eq:exrwg}
		&e_{x,n}^{\rm TE} = \hspace*{3mm}\frac{B_{pq}}{k_{pq}}k_{yq}\cos(k_{xp}\tilde{x})\sin(k_{yq}\tilde{y}), \\
		&e_{y,n}^{\rm TE} = -\frac{B_{pq}}{k_{pq}}k_{xp}\sin(k_{xp}\tilde{x})\cos(k_{yq}\tilde{y}),
		\end{align}
		while, if the $n$-th mode is TM, it follows 
		\begin{align}
		&e_{x,n}^{\rm TM} = \frac{B_{pq}}{k_{pq}}k_{xp}\cos(k_{xp}\tilde{x})\sin(k_{yq}\tilde{y}), \\
		&e_{y,n}^{\rm TM} = \frac{B_{pq}}{k_{pq}}k_{yq}\sin(k_{xp}\tilde{x})\cos(k_{yq}\tilde{y}).
		\end{align}
		Moreover, for TM modes, we have also to define
		\begin{equation}\label{eq:ezrwg}
		e_{z,m}^{\rm TM} = B_{pq}\sin(k_{xp}\tilde{x})\sin(k_{yq}\tilde{y}),
		\end{equation}
		for $m=1,\ldots, N_{{\rm TM}}$.
		In \eqref{eq:exrwg}-\eqref{eq:ezrwg}, we have set
		\begin{align}
		\tilde{x} = x - a_0/2,& \hspace{2em}\tilde{y} = y - b_0/2, \\
		k_{xp}=\frac{p\pi}{a_0},&\hspace{2em}k_{yq}=\frac{q\pi}{b_0}, \\
		k_{pq} = \sqrt{k_{xp}^2+k_{yq}^2},&\hspace{2em}
		B_{pq} = \sqrt{\frac{\varepsilon_{p0}\varepsilon_{q0}}{a_0b_0}},
		\end{align}
		where $\varepsilon$ is a function of $p$ or $q$, being
		\begin{align}
		&\varepsilon_{p0} = \left\{
		\begin{matrix}
		1\hspace{1em}{\textrm{for }} p=0  \\
		2\hspace{1em}{\textrm{for }} p> 0, & 
		\end{matrix}
		\right. \\[2mm]
		&\varepsilon_{q0} = \left\{
		\begin{matrix}
		1\hspace{1em}{\textrm{for }} q=0  \\
		2\hspace{1em}{\textrm{for }} q> 0. & 
		\end{matrix}
		\right.
		\end{align}
		The $N_{\rm M}$ modes are selected as the first $N_{\rm M}$ modes in the increasing sequence of the associated eigenvalue $k_{pq}$.

	
	\null
	
	\bibliographystyle{IEEEtran}
	\bibliography{IEEEabrv,./paper.bib}
	
	
\end{document}